\documentclass[11pt,a4paper]{article}

\usepackage[margin=2.5cm]{geometry}
\usepackage[T1]{fontenc}
\usepackage[utf8]{inputenc}
\usepackage{lmodern}
\usepackage{amsmath,amssymb,amsthm}
\usepackage{graphicx}
\usepackage{booktabs}
\usepackage{microtype}
\usepackage{siunitx}
\usepackage[numbers]{natbib}
\usepackage{hyperref} 
\usepackage{float}
\usepackage[utf8]{inputenc}
\DeclareUnicodeCharacter{2032}{\ensuremath{\prime}}

\hypersetup{
    colorlinks=true,
    linkcolor=blue,
    citecolor=blue,
    urlcolor=blue
}

\title{Information funnels and multiscale gap-space dynamics in Kaprekar's routine}

\author{Christoph D.~Dahl\\[2mm]
\small Graduate Institute of Mind, Brain and Consciousness,\\[-1mm]
\small Taipei Medical University, Taipei, Taiwan\\[1mm]
\small \texttt{christoph.d.dahl@gmail.com}}
\date{}

\begin{document}
\maketitle

\begin{abstract}
Kaprekar's routine, i.e., sorting the digits of an integer in ascending and descending order and subtracting the two, defines a finite deterministic map on the state space of fixed-length digit strings. While its attractors (such as $495$ for $D=3$ and $6174$ for $D=4$) are classical, the global information-theoretic structure of t§ve analysis is carried out for $D\in\{3,4,5,6\}$. For each $D$, all states are enumerated, their attractors and convergence distances are obtained, and the induced distribution over attractors across iterations is used to construct ``entropy funnels''. Despite the combinatorial growth of the state space, average distances remain small and entropy decays rapidly before entering a slow tail. Permutation symmetry is then exploited by grouping states into digit multisets and, in a further reduction, into low-dimensional digit-gap features. On this gap space, Kaprekar's routine induces a first-order Markov approximation whose transition structure, stationary distribution and drift fields are characterised, showing that simple gap features strongly constrain the dynamics for $D=3$ but lose predictive power as $D$ increases.
\end{abstract}

\section*{Introduction}
Iterated digit transforms provide simple yet surprisingly rich examples of finite dynamical systems. Among them, Kaprekar's routine occupies a special place: starting from a $D$-digit integer with at least two distinct digits, its digits are sorted into descending and ascending order, interpreted as integers, and subtracted. Formally, for a $D$-digit state $x$ with digits $(d_1,\dots,d_D)$, define
\begin{equation}
\label{eq:desc_asc}
\begin{aligned}
\mathrm{desc}(x) &= \text{integer formed by sorting the digits of $x$ in descending order},\\
\mathrm{asc}(x)  &= \text{integer formed by sorting the digits of $x$ in ascending order}.
\end{aligned}
\end{equation}
The Kaprekar map is then given by
\begin{equation}
    K(x) = \mathrm{desc}(x) - \mathrm{asc}(x),
    \label{eq:kaprekar_map}
\end{equation}
iterated on the finite state space of $D$-digit integers (with leading zeros allowed). For $D=4$, this process famously converges to the attractor $6174$ for almost all initial conditions, while $495$ plays an analogous role for $D=3$ \citep{kaprekar1949another, kaprekar1955interesting, nishiyama2006mysterious}. These facts are well known in recreational mathematics, and classical work has established existence, uniqueness and basic properties of Kaprekar attractors in various bases and digit lengths \citep{trigg1971kaprekar, trigg1972kaprekar, eldridge1988determination, prichett1981determination, walden2005searching, dolan2011classification}. Beyond these combinatorial results, however, relatively little is known about the \emph{global} organisation of Kaprekar dynamics, its information-theoretic signatures, or how these properties depend on the digit length~$D$. Seen as a finite dynamical system, Kaprekar’s routine raises a few concrete questions: (1) How are basins of attraction distributed as $D$ increases, and how dominant is the largest basin? (2) How quickly does uncertainty about the eventual attractor collapse under iteration, when starting from a uniform prior over states? (3) Are there simple low-dimensional features of the digits that control whether a state is ``easy'' or ``hard'' to reach? (4) Can the dynamics be captured, at least approximately, by a stochastic process on a coarse-grained state space?
Recent mathematical work has addressed structural and asymptotic properties of Kaprekar-type maps in a variety of settings, including $b$-adic generalisations, bounds on Kaprekar constants, and detailed analyses of two- and three-digit routines \cite{yamagami20182,wang20212,nuez2021symmetries,devlin2020maximum}. These studies highlight the richness of the underlying number-theoretic structure, but they do not aim to characterise the induced dynamics in information-theoretic or probabilistic terms.\\
~\\
Here a different angle is taken: Instead of trying to obtain closed-form expressions for constants or loops, the Kaprekar map is treated as a finite directed graph, and questions are posed about its global statistics and coarse-grained descriptions: how entropy contracts, how permutation symmetries can be exploited, and how far simple `gap' features go in explaining the dynamics.
Standard tools from information theory and Markov-chain analysis \cite{shannon1948mathematical,cover1999elements,norris1998markov} are used. In the present work these questions are addressed for $D\in\{3,4,5,6\}$ by combining exhaustive enumeration with information-theoretic and statistical tools. The map $K$ in~\eqref{eq:kaprekar_map} is treated as a deterministic update on a finite directed graph, and its structure is analysed at three complementary levels. First, at the level of individual states and attractors, all $D$-digit states with at least two distinct digits are enumerated, their attractors and distances (numbers of iterations) to convergence are obtained, and basin sizes are quantified. This representation makes it possible to construct \emph{entropy funnels}: starting from a uniform distribution over states, the evolving distribution over attractors among those trajectories that have already converged is followed across iterations and the corresponding Shannon entropy is measured, providing a quantitative notion of uncertainty reduction under the dynamics. Second, the permutation symmetry of digit strings is exploited by grouping states into equivalence classes with identical digit multisets. This reduction yields a more compact representation that still respects the combinatorial constraints of the map. Multiset-level statistics are used to characterise how class sizes are distributed, how far typical classes lie from their attractors, and how different attractors are assembled from contributions of many small versus few large classes. Third, a low-dimensional description is introduced in terms of simple digit-gap features,
\begin{equation}
    g_1 = \max(d) - \min(d), \qquad
    g_2 = d_{(2)} - d_{(3)},
    \label{eq:gap_features}
\end{equation}
where $d_{(k)}$ denotes the $k$-th largest digit. These features define a discrete ``gap space''. From the exact deterministic dynamics one can estimate empirical one-step transition probabilities between gap states, which define a first-order Markov chain on this coarse-grained space: each step maps the empirical distribution over gap states to a new one. On this gap space, transition probabilities, stationary distributions, flow fields and drift statistics are estimated, and the resulting structure is related back to the basins and distances in the original state space. In addition, a simple linear regression framework is used to probe how far gap-based and aggregate digit features can account for the distance-to-attractor, viewed as a notion of ``difficulty'' of reaching a given state.
These three levels of description (states, multisets, and gap space) give a multi-scale view of Kaprekar’s routine: from individual trajectories and basins to symmetry-reduced classes and finally to a very low-dimensional Markov picture. The same approach should extend to other digit-based transforms and, more generally, to deterministic maps where one cares about coarse-grained information flow.

\section*{Methods}
Fix a digit length $D \ge 3$ and work in base~10.
Let $\mathcal{S}_D$ denote the set of $D$-digit states with at least two distinct digits, allowing leading zeros.
Write $x\in\mathcal{S}_D$ as a $D$-tuple of digits $(d_1,\dots,d_D)$.
The descending and ascending orderings $\mathrm{desc}(x)$ and $\mathrm{asc}(x)$ are defined as in~\eqref{eq:desc_asc}, and the Kaprekar map is given by~\eqref{eq:kaprekar_map}; when restricted to digit length $D$ it is written as $K_D : \mathcal{S}_D \to \mathcal{S}_D$.
The map $K_D$ is iterated and $K_D^t$ denotes the $t$-fold composition.
A state $a\in\mathcal{S}_D$ is an \emph{attractor} if $K_D(a)=a$, and its (forward) \emph{basin} is
\begin{equation}
\mathcal{B}(a) := \{x\in\mathcal{S}_D : K_D^t(x)=a \text{ for some } t\ge 0\}.
\label{eq:basin}
\end{equation}
For each $x\in\mathcal{S}_D$ there exists $t\ge 0$ such that $K_D^t(x)$ is eventually periodic.
In the computations reported here all attracting cycles encountered are fixed points, i.e.\ $K_D$-attractors in the above sense.
For $x\in\mathcal{B}(a)$ the \emph{distance-to-attractor} is defined as
\begin{equation}
\mathrm{dist}(x) := \min\{t\ge 0 : K_D^t(x)=a\}.
\label{eq:dist}
\end{equation}
It is convenient to write $\alpha(x)$ for the attractor reached from $x$ and $\tau(x):=\mathrm{dist}(x)$ for the corresponding convergence time.
For each $D\in\{3,4,5,6\}$ the set $\mathcal{S}_D$ is enumerated, all attractors are identified, basin sizes $|\mathcal{B}(a)|$ are computed, and $\mathrm{dist}(x)$ is recorded for all states.\\[1ex]
To quantify information funnels an initial distribution that is uniform over all non-trivial states for a given~$D$ is considered.
For each iteration $t\ge 0$ define the subset of states that have already reached an attractor by time~$t$,
\begin{equation}
\mathcal{A}_t := \{x\in\mathcal{S}_D : \tau(x) \le t\}, \qquad N_t := |\mathcal{A}_t|.
\label{eq:At}
\end{equation}
Among these converged states, the empirical distribution over attractors at time $t$ is
\begin{equation}
p_t(a) = \frac{1}{N_t}\bigl|\{x\in\mathcal{A}_t : \alpha(x)=a\}\bigr|,
\label{eq:pt}
\end{equation}
where the sum over $a$ runs over all attractors for the given $D$.
The Shannon entropy
\begin{equation}
H_t = -\sum_a p_t(a)\log_2 p_t(a)
\label{eq:entropy}
\end{equation}
is then computed as a function of iteration $t$.
For small $t$ the distribution $p_t$ is dominated by attractors that are reached quickly; as $t$ increases and more trajectories converge, $p_t$ approaches the basin-size distribution.
Plotting $H_t$ yields a raw entropy funnel; normalising by $H_{t^\ast}$, where $t^\ast$ is chosen large enough that all states have converged for the given $D$, gives the relative reduction in uncertainty relative to the final attractor distribution.\\
Because permuting the digits of $x$ does not affect the outcome of a Kaprekar step, many states form equivalence classes with identical digit multisets.
Formally, an equivalence relation is defined on $\mathcal{S}_D$ by $x\sim y$ if their digits coincide as multisets.
The equivalence classes are in bijection with digit multisets, and are referred to as \emph{multiset classes}.
Each class has a \emph{size} (number of distinct permutations) and a well-defined mean distance-to-attractor obtained by averaging $\mathrm{dist}(x)$ over the states in the class.
For each $D$ all digit multisets, their class sizes, and their mean distances are enumerated.
It is also recorded, for each attractor, how many states in its basin arise from each multiset class.\\
~\\
For each state $x\in\mathcal{S}_D$ simple digit features are computed.
In the present work the focus is on the gap features introduced in~\eqref{eq:gap_features}, written explicitly as
\begin{equation}
g_1(x) = \max(d) - \min(d), \qquad
g_2(x) = d_{(2)} - d_{(3)},
\label{eq:gap_features_methods}
\end{equation}
where $d_{(1)}\ge d_{(2)}\ge \dots\ge d_{(D)}$ are the sorted digits of $x$.
These define a discrete set
\begin{equation}
\mathcal{G}_D := \bigl\{(g_1(x),g_2(x)) : x\in\mathcal{S}_D\bigr\}
\label{eq:gap_space}
\end{equation}
of possible gap pairs $(g_1,g_2)$.\\
Kaprekar's routine induces a stochastic map on this grid:
for each gap state $g=(g_1,g_2)\in\mathcal{G}_D$, consider the collection of underlying states
\begin{equation}
\mathcal{S}_D(g) := \{x\in\mathcal{S}_D : (g_1(x),g_2(x)) = g\}.
\label{eq:state_class}
\end{equation}
One Kaprekar step is applied to each $x\in\mathcal{S}_D(g)$, the successor gaps $g'(x)=(g_1(K_D(x)),g_2(K_D(x)))$ are recorded, and empirical transition frequencies
\begin{equation}
P_D(g\to h) := \frac{1}{|\mathcal{S}_D(g)|}\bigl|\{x\in\mathcal{S}_D(g): g'(x)=h\}\bigr|
\label{eq:transition}
\end{equation}
are estimated.
These empirical probabilities define a first-order Markov chain on $\mathcal{G}_D$ with transition matrix $P_D$, which provides a coarse-grained approximation to the projected dynamics on gap space (the true projected process need not be strictly Markov).
The following quantities are computed: (1) the stationary distribution $\pi_D$ satisfying $\pi_D P_D = \pi_D$; (2) the empirical distribution of gap states under the uniform prior on $\mathcal{S}_D$; (3) average changes $\Delta g_1$ and $\Delta g_2$ per step for each gap state.\\
To link local digit structure to distance-to-attractor a simple regression analysis is performed.
For each $D$ up to $N = 50\,000$ states (or all states when fewer are available) are sampled from $\mathcal{S}_D$ and, for each state, $g_1$ and $g_2$ as defined above, the digit sum $\sum_i d_i$, and the variance of digits $\mathrm{Var}(d)$ are computed.
All features are standardised (zero mean, unit variance).
A linear model
\begin{equation}
  \widehat{\mathrm{dist}}(x)
    = \beta_0
    + \beta_1\, g_1(x)
    + \beta_2\, g_2(x)
    + \beta_3 \sum_{i=1}^D d_i
    + \beta_4 \operatorname{Var}(d).
\label{eq:regression}
\end{equation}
is then fit by least squares, and, for each $D$, the coefficient of determination $R^2$, the root mean squared error (RMSE), and the learned weights~$\beta$ are reported. For the comparisons between ``easy'' and ``hard'' states in Figure~\ref{fig:easy_vs_hard}, states are ranked by $\mathrm{dist}(x)$ and the fastest and slowest deciles (smallest and largest $10\%$ of distances) are selected. Digit features are then averaged separately for these two groups.
\\
All computations were performed in MATLAB (Mathworks\textsuperscript{\tiny\textregistered}, Natick, MA, USA), using integer-valued operations for the Kaprekar map that are exactly represented in double-precision arithmetic for the ranges of $D$ considered, and double-precision arithmetic for derived quantities such as entropies.
Enumeration of $\mathcal{S}_D$ and construction of the gap-space Markov chain are exact for each $D\in\{3,4,5,6\}$.
Where sampling is used (in the regression analysis), the corresponding sample size is stated.

\section*{Results}
\begin{figure}[t]
    \centering
    \includegraphics[width=112mm]{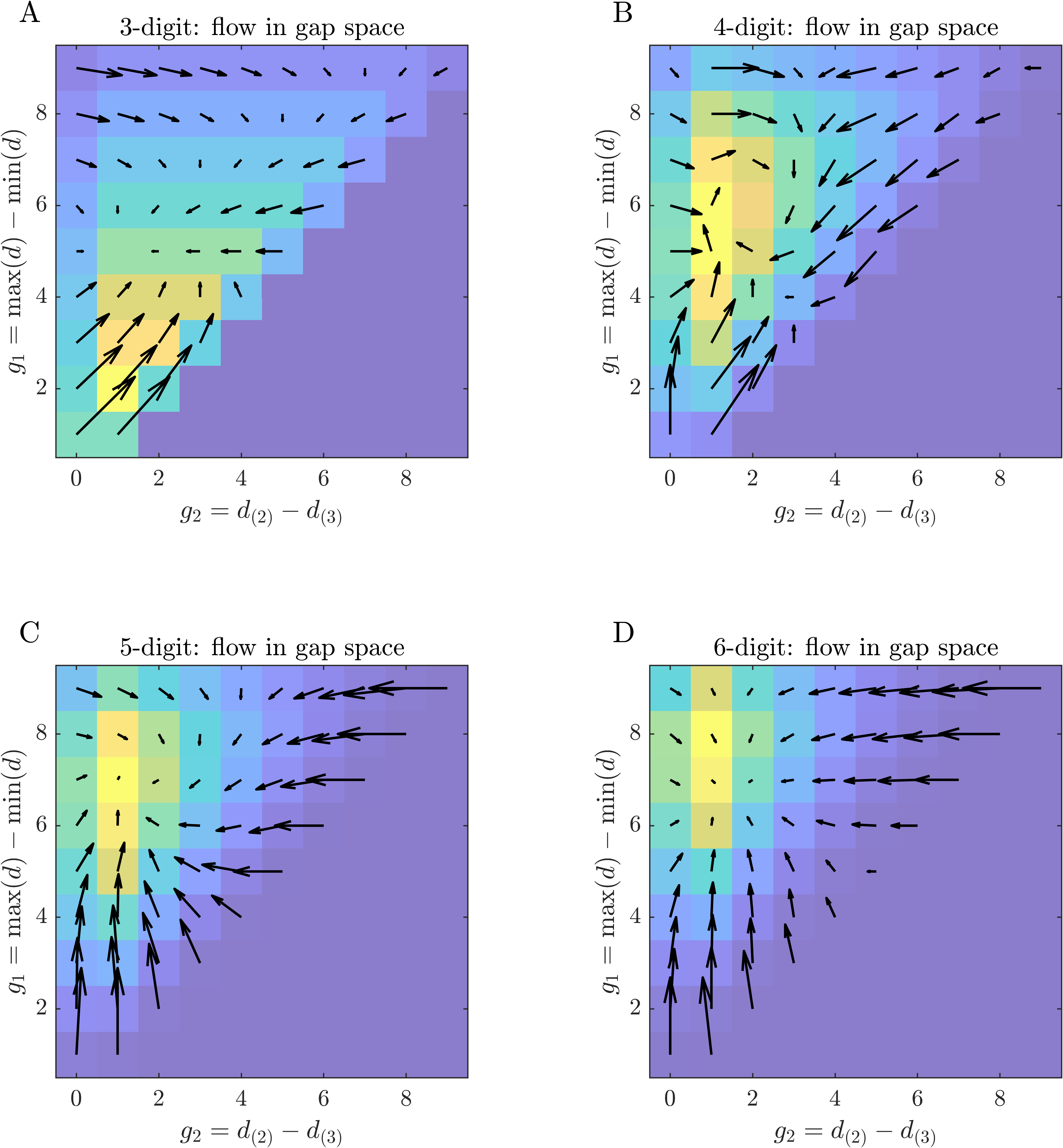}
    \caption{\textbf{Flow fields in gap space for $D=3$--$6$.} Panels (\textbf{A})–(\textbf{D}) correspond to $D=3,4,5,6$, respectively. Each panel shows the occupancy of gap states $(g_1,g_2)$ under a uniform prior
over $\mathcal{S}_D$ (colour scale) together with the corresponding average one-step drift induced by a single Kaprekar step. Arrows indicate the mean change in gap coordinates per iteration, with horizontal and vertical
components corresponding to $\Delta g_2$ and $\Delta g_1$, respectively.}
    \label{fig:gap_flow_allDs}
\end{figure}

\begin{figure}[t]
    \centering
    \includegraphics[width=80mm]{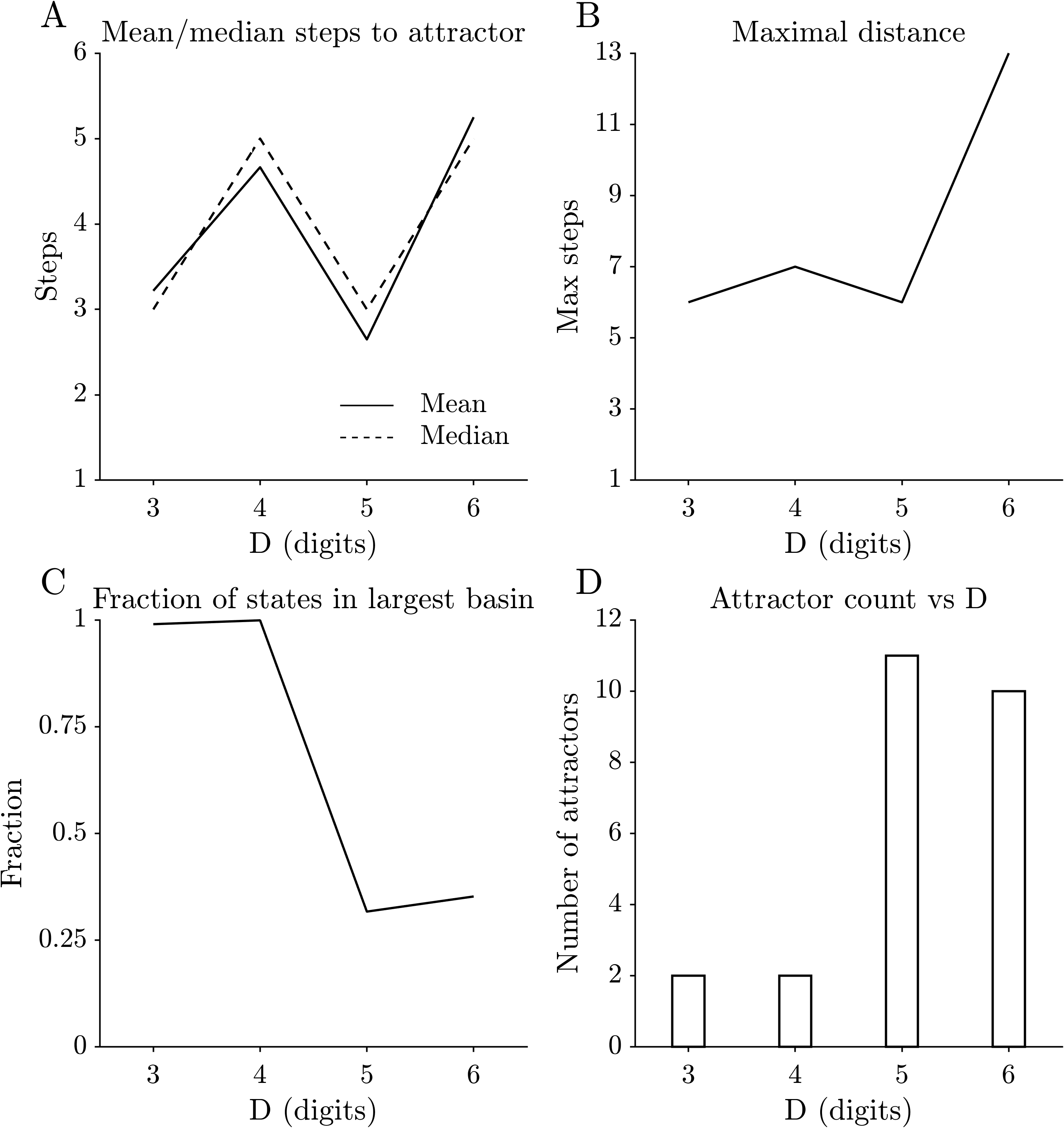}
     \caption{\textbf{Global structure of Kaprekar dynamics as a function of digit length $D$.}
(\textbf{A}) Mean (solid) and median (dashed) number of iterations required to reach an attractor.
(\textbf{B}) Maximal distance over all starting states.
(\textbf{C}) Fraction of the state space contained in the largest basin of attraction.
(\textbf{D}) Number of distinct attractors.
}
    \label{fig:global_landscape}
\end{figure}

\begin{figure}[t]
    \centering
	\includegraphics[width=80mm]{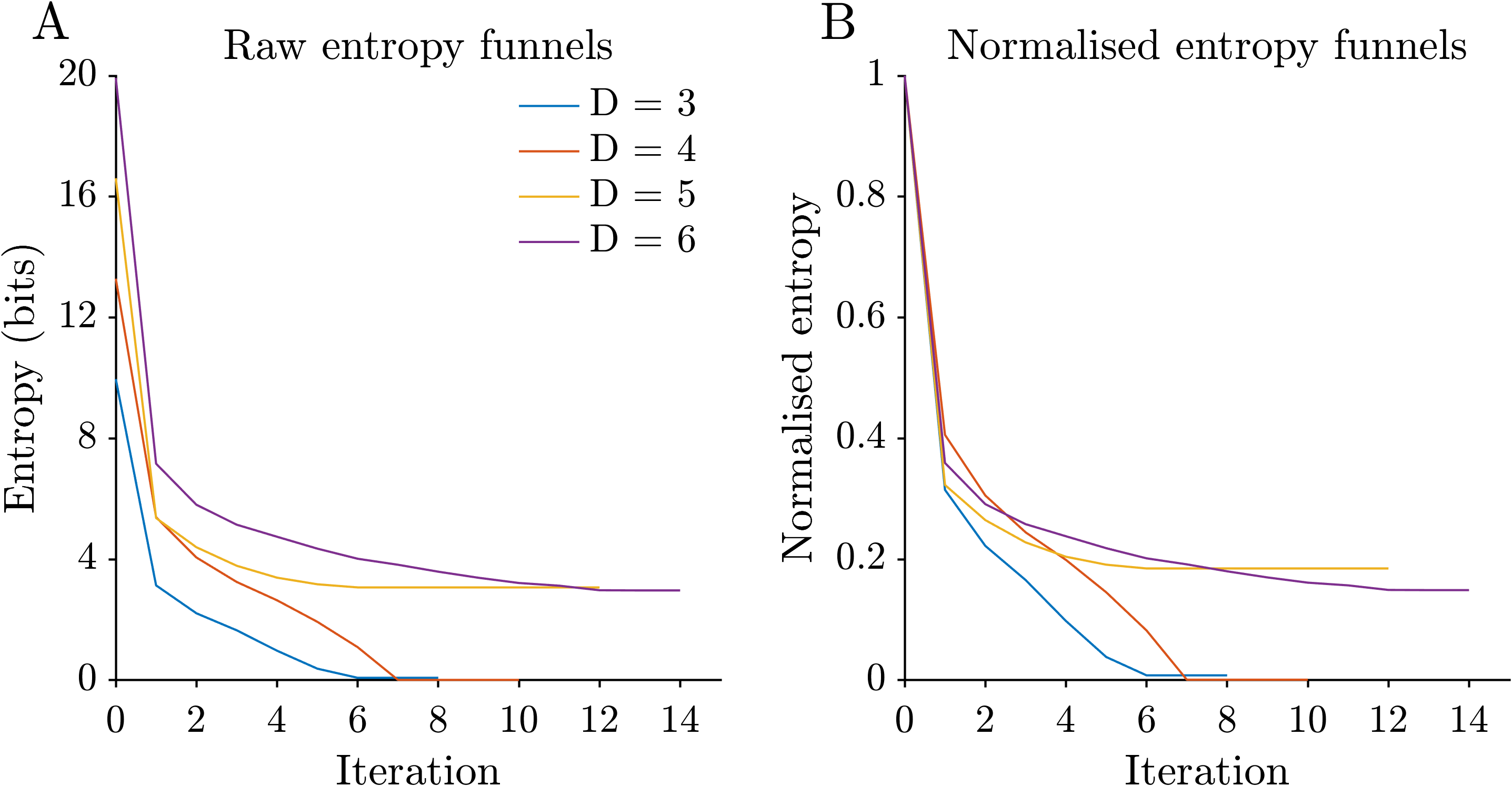}
    \caption{\textbf{Entropy funnels across digit lengths.}
(\textbf{A}) Raw Shannon entropy of the induced distribution over attractors, conditional on trajectories that have converged by iteration $t$, for each $D$.
(\textbf{B}) The same curves normalised by their final values $H_{t^\ast}$ (with $t^\ast$ chosen large enough that all states have converged), showing relative uncertainty about the final attractor.
All digit lengths exhibit a rapid early collapse of uncertainty, followed by a slower tail governed by the distribution of basin sizes.}
 \label{fig:entropy_funnels}
\end{figure}

\begin{figure}[t]
    \centering
	\includegraphics[width=80mm]{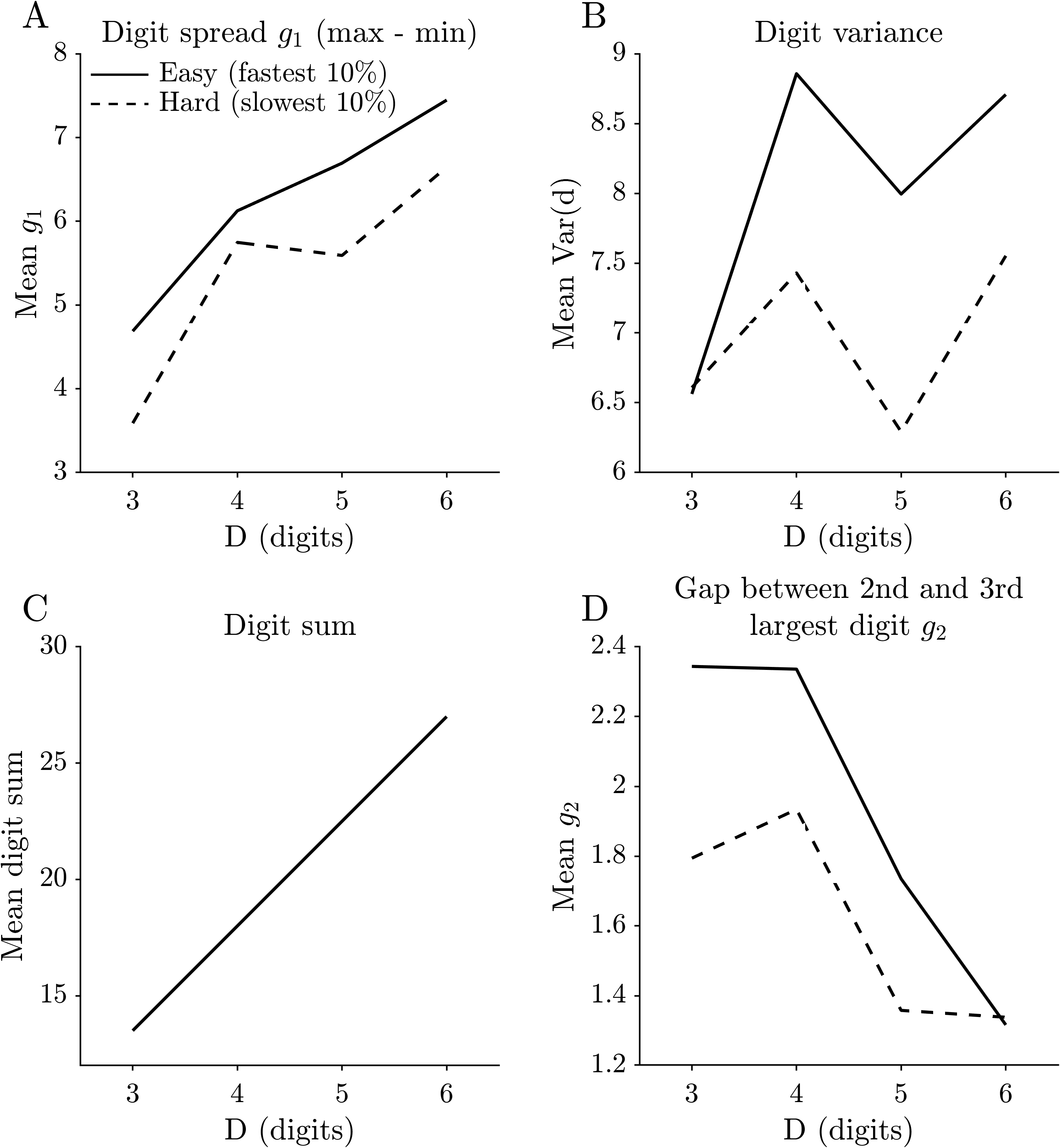}
    \caption{\textbf{Digit features for easy and hard Kaprekar states.}
For each digit length $D$, states are split into the fastest and slowest deciles of distance-to-attractor.
Panels (\textbf{A})–(\textbf{D}) show, respectively, the mean digit spread $g_1 = \max(d) - \min(d)$, digit variance $\mathrm{Var}(d)$, digit sum, and the gap between the second and third largest digits $g_2$ for the two groups.
Easy states consistently exhibit larger $g_1$ and, for some $D$, higher digit variance, whereas differences in digit sum are largely driven by the change in $D$.}
    \label{fig:easy_vs_hard}
\end{figure}

\begin{figure}[t]
    \centering
	\includegraphics[width=80mm]{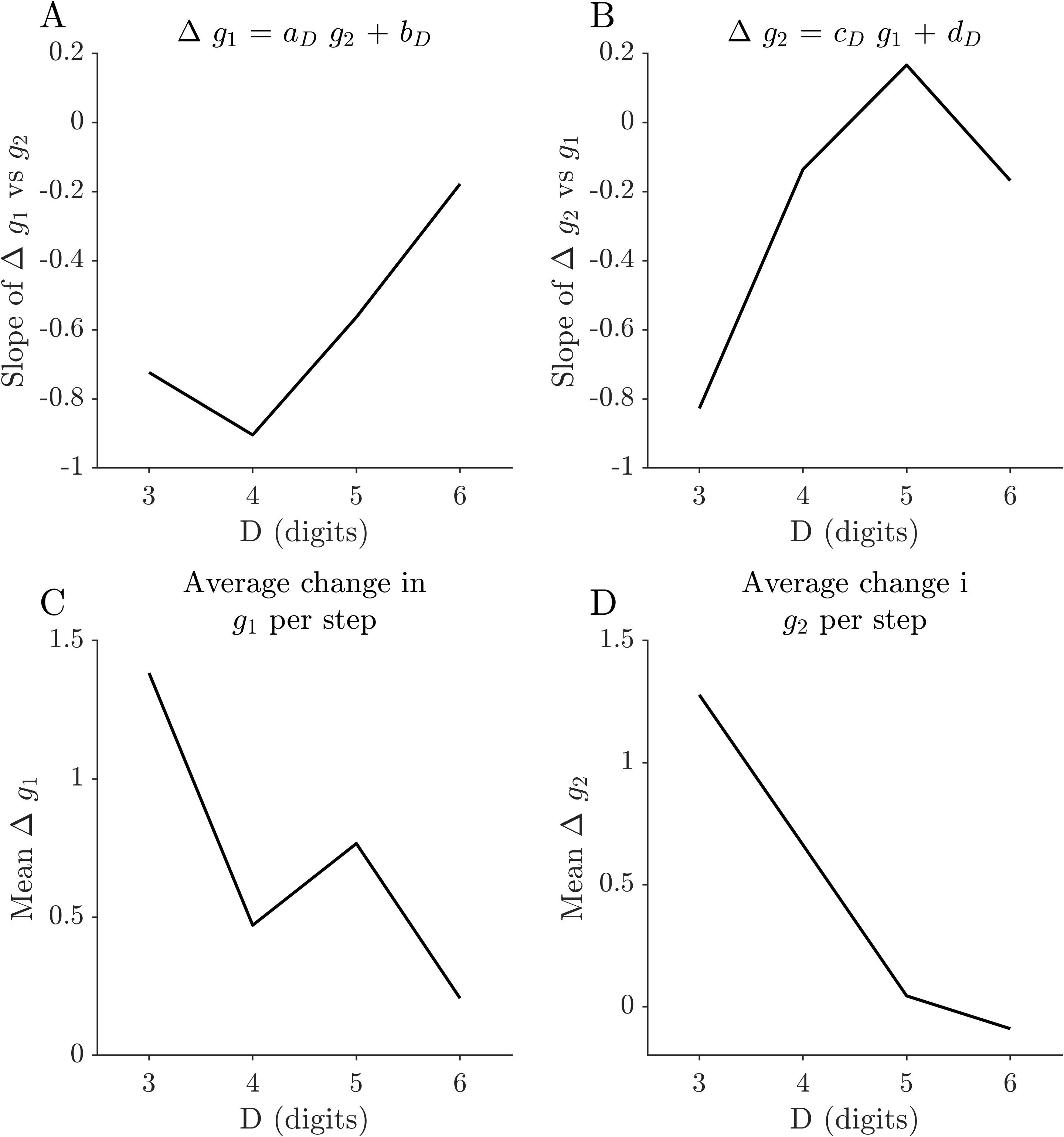}
    \caption{\textbf{Cross-digit summary of average gap drift.} Panels (\textbf{A}) and (\textbf{B}) show the slopes of the linear relations $\Delta g_1 \approx a_D g_2 + b_D$ and $\Delta g_2 \approx c_D g_1 + d_D$, respectively, as a function of $D$, obtained by ordinary least-squares fits across all occupied gap states. The negative slopes indicate a robust coupling between the two gap coordinates that weakens with increasing $D$. Panels (\textbf{C}) and (\textbf{D}) show the average drifts $\langle\Delta g_1\rangle$ and $\langle\Delta g_2\rangle$ per step, revealing a systematic tendency to increase $g_1$ and decrease $g_2$ that becomes smaller in magnitude for higher $D$.}
    \label{fig:gap_drift}
\end{figure}

\begin{figure}[t]
    \centering
	\includegraphics[width=43mm]{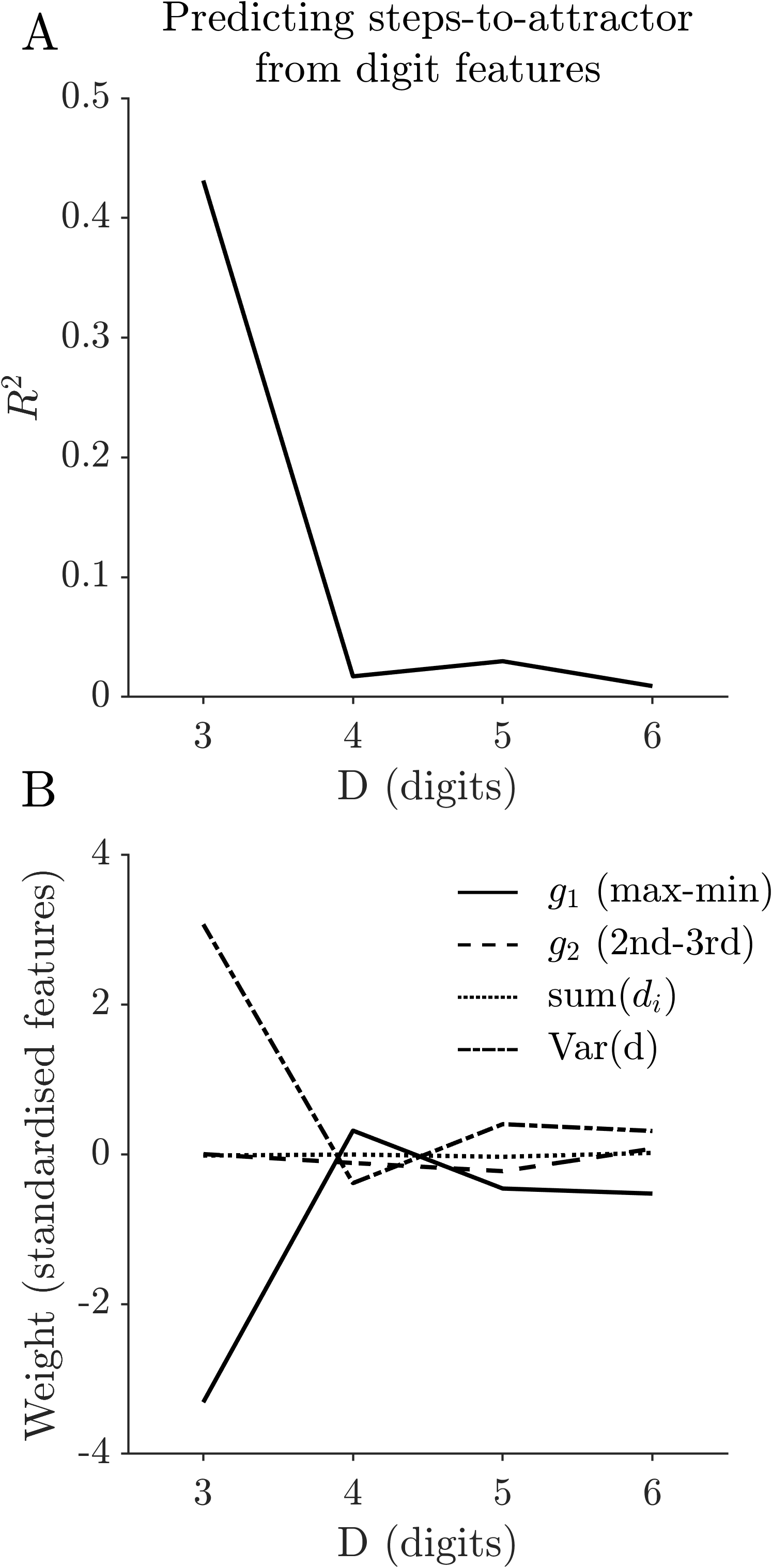}
    \caption{\textbf{Predicting distance-to-attractor from simple digit features.}
Panel (\textbf{A}) shows the coefficient of determination $R^2$ for linear regression models predicting the number of Kaprekar steps from standardised features $(g_1,g_2,\text{digit sum},\text{digit variance})$ as a function of digit length~$D$. Predictive power is substantial for $D=3$ but drops sharply for $D \ge 4$ in this simple linear model.
Panel (\textbf{B}) shows the learned regression weights for each feature across $D$, indicating that large digit spread $g_1$ is strongly associated with faster convergence only for three digits, while for larger $D$ the relationship between these local features and global distance-to-attractor becomes weak within this linear framework.}
    \label{fig:predictability}
\end{figure}

Before scalar summary statistics are introduced, the coarse-grained flow of Kaprekar dynamics in gap space is visualised. Figure~\ref{fig:gap_flow_allDs} shows, for each digit length $D \in \{3,4,5,6\}$, the occupancy of gap states $(g_1,g_2)$ under a uniform prior over $\mathcal{S}_D$ together with the corresponding average one-step drift vectors. Colours encode how many states realise a given gap pair, and arrows indicate the mean change $(\Delta g_2,\Delta g_1)$ produced by a single Kaprekar step. For three digits, trajectories concentrate along a narrow band and converge towards a single high-occupancy region. For larger $D$ the occupied region of gap space broadens and the flow becomes more heterogeneous, with weaker and more dispersed drift, foreshadowing the more fragmented basin structure quantified below.
Figure~\ref{fig:global_landscape} summarises how the global structure of Kaprekar dynamics changes with the number of digits~$D$. For each $D$ the following summary quantities are reported: the number of distinct attractors, the fraction of states lying in the largest basin of attraction, the mean and median distance-to-attractor, and the maximal distance observed. Here the distance-to-attractor measures how many iterations of the Kaprekar map are needed before the trajectory reaches its final fixed point. Despite the combinatorial explosion of the state space as $D$ increases, the average number of iterations needed to reach an attractor remains small for all values of $D$ considered. By contrast, the maximal distance, the dominance of the largest basin, and the number of distinct attractors all vary systematically with~$D$. For $D=3$ and $D=4$ a single attractor dominates the dynamics, in the sense that most initial states eventually flow into one large basin.
For $D=5$ and $D=6$ the picture is more fragmented: the largest basin occupies a much smaller fraction of the state space and many additional, smaller attractors appear.
\\
~\\
Information contraction under Kaprekar iteration can be summarised by \emph{entropy funnels}, which track how uncertainty about the eventual attractor decreases over time. For each digit length $D$, Figure~\ref{fig:entropy_funnels} shows the Shannon entropy of the induced distribution over attractors, among those trajectories that have converged by iteration $t$, as a function of iteration, together with a normalised version. High entropy corresponds to a situation in which many attractors are still plausible, whereas low entropy indicates that almost all initial states have effectively committed to a small subset of attractors. Raw entropy (in bits) decays rapidly in the first few iterations and then flattens as trajectories approach their attractors and the basin-size distribution is revealed. The normalised curves show that, for all $D$, the bulk of uncertainty about the eventual attractor is resolved within roughly five iterations, but the final residual entropy depends strongly on the number and relative sizes of basins. For $D=3$ and $D=4$ the normalised entropy drops close to zero, reflecting the near-complete dominance of a single attractor. For $D=5$ and $D=6$ the decay is less complete, consistent with the proliferation of attractors and a more even distribution of basin sizes, so that some uncertainty about the final attractor remains even after many iterations.\\
~\\
Factoring states into digit multisets provides a compact summary of the combinatorial structure of the map (Figure~\ref{fig:multiset_structure}). A digit multiset records which digits appear and how often, but ignores their order, so all permutations of the same digits belong to the same multiset class. This perspective separates purely combinatorial effects (how many permutations a given multiset admits) from dynamical effects (how quickly states built from that multiset converge). For each digit length $D\in\{3,4,5,6\}$, the left-hand panel shows the distribution of multiset sizes, i.e.\ the number of distinct permutations in each class, both as absolute counts (solid line, left axis) and as normalised probabilities (dashed line, right axis). As $D$ increases, the number of multiset classes grows and the distributions shift towards larger sizes with heavier upper tails, indicating the appearance of rare digit patterns with many distinct permutations. The right-hand panels display the corresponding distributions of mean distance-to-attractor per multiset, obtained by averaging $\mathrm{dist}(x)$ over all states in a class. For all $D$, most multisets are associated with short mean distances, but the tails extend further for $D=5$ and $D=6$, reflecting a growing minority of digit patterns that are systematically linked to longer transients.
This multiset-level view shows that part of the complexity at higher digit lengths arises from an increasingly uneven allocation of basin structure across combinatorial classes.
\begin{figure}[t]
    \centering
	\includegraphics[width=96mm]{figures/Fig3_multiset_structure.png}
    \caption{\textbf{Multiset structure of Kaprekar dynamics across digit lengths.}
Panels (\textbf{A}), (\textbf{C}), (\textbf{E}), and (\textbf{G}) show the distributions of multiset sizes, i.e.\ the number of distinct permutations in each digit-multiset class, for $D = 3,4,5,6$, respectively.
Panels (\textbf{B}), (\textbf{D}), (\textbf{F}), and (\textbf{H}) show the corresponding distributions of mean distance-to-attractor per multiset class (average number of Kaprekar iterations taken by all states sharing the same multiset of digits).
Solid black lines show absolute counts (left $y$-axis); dashed red lines show the corresponding normalised probabilities (right $y$-axis).
As $D$ increases, both the number and typical size of multiset classes grow sharply, and the distributions of multiset-averaged distances broaden and develop longer tails, indicating increasingly heterogeneous convergence behaviour. Note that most classes have small mean distance (fewer than three steps), with a long right tail for $D=5$ and $D=6$.}
	 \label{fig:multiset_structure}
\end{figure}
For each digit length $D$, states are ranked by their distance-to-attractor and divided into two extreme groups: the fastest $10\%$ are referred to as ``easy'' states and the slowest $10\%$ as ``hard'' states. Comparing their digit features reveals systematic differences (Figure~\ref{fig:easy_vs_hard}). Easy states tend to have a larger overall digit spread $g_1 = \max(d) - \min(d)$, and for some $D$ they also show higher digit variance, indicating that configurations whose digits are more widely dispersed are typically closer to an attractor in the dynamical sense. Differences in digit sum are much less informative: because the total sum of digits naturally increases with $D$, shifts in this feature are largely driven by digit length rather than by dynamical difficulty and should therefore be interpreted with caution.\\
~\\
In gap space, Kaprekar’s routine induces a Markov chain with highly structured flow fields (Figure~\ref{fig:gap_flow_allDs}). For $D=3$, the occupied region of gap space forms a narrow wedge constrained by $g_2 \le g_1$, with the highest occupancy at moderate $g_1$ and small-to-moderate $g_2$. The corresponding drift vectors show a coherent flow towards larger $g_1$ and moderate $g_2$, foreshadowing the basin structure quantified below. For larger $D$, the flow patterns become progressively more asymmetric, and the stationary distributions place most of their mass in regions with large $g_1$ and moderate $g_2$, indicating a bias towards states with a wide overall digit spread but only a moderate gap between the second and third largest digits. Average drift vectors show a robust tendency to increase $g_1$ and decrease $g_2$, with drift magnitudes decreasing as $D$ grows. To summarise these trends across digit lengths, one-step changes $\Delta g_1$ are regressed on $g_2$ and $\Delta g_2$ on $g_1$ for each $D$, and mean changes per step are computed (Figure~\ref{fig:gap_drift}). The slopes of the linear relations $\Delta g_1 \approx a_D g_2 + b_D$ and $\Delta g_2 \approx c_D g_1 + d_D$ are negative for all $D$, but their magnitude decreases with $D$, indicating that the coupling between the two gap coordinates weakens in higher-digit systems.
From a more intuitive perspective, the dynamics tend to push states towards regions of gap space where one large overall digit range coexists with a more balanced configuration among the middle digits, but this directional bias becomes less pronounced as $D$ grows.
The regression analysis (Figure~\ref{fig:predictability}) quantifies how well this particular choice of simple digit features predicts distance-to-attractor. For $D=3$, the linear model explains approximately $R^2 \approx 0.43$ of the variance in distance with a root mean squared error (RMSE) of about $1.08$ steps over $N=998$ states. For $D=4,5,6$, the $R^2$ values drop to roughly $0.02$--$0.03$ with markedly larger RMSE, indicating that, beyond three digits, only a very small fraction of the variability in distance-to-attractor can be accounted for by these four standardised features within a linear model. In other words, for $D\ge 4$, the link between this particular low-dimensional linear summary of local digit structure and global convergence time is weak. The learned regression weights highlight how the influence of individual features changes with $D$. For $D=3$, a large digit spread $g_1$ is strongly associated with faster convergence, and digit variance carries a clear positive weight, suggesting that states with more heterogeneous digits tend to approach their attractors quickly. For larger $D$, the weights shrink in magnitude and fluctuate in sign across features, consistent with a more entangled and higher-dimensional dependence of distance-to-attractor on the underlying digit configuration. From a broader perspective, these findings show that while simple low-dimensional linear descriptors are informative for three-digit Kaprekar dynamics, they capture only a small part of the variability for $D\ge 4$; richer feature sets or nonlinear models might recover additional structure but are beyond the scope of the present analysis.\\

\section*{Discussion}
The results above describe Kaprekar's routine, for $D = 3$--$6$, as a finite dynamical system seen through several information-theoretic lenses. A very simple digit transform already produces a surprisingly layered structure:
shallow typical transients but long tails, a shift from one dominant basin to many smaller ones, and a useful, but ultimately limited, low-dimensional gap description.
At the coarsest level, the global summaries in Figure~\ref{fig:global_landscape} indicate that Kaprekar dynamics remains shallow on average across all digit lengths considered: typical trajectories reach an attractor in only a few iterations, despite the combinatorial growth of the state space. This shallow behaviour coexists with marked changes in extremal and structural quantities. The maximal distance-to-attractor increases with $D$, the dominance of the largest basin decreases, and the number of distinct attractors rises. For three and four digits, the dynamics are effectively governed by a single large basin, whereas for five and six digits the state space splits into many smaller basins with more heterogeneous sizes. This transition from ``one big funnel'' to a more fragmented landscape is a first indication that, within the range $D=3$--$6$ examined here, the global organisation of the map shifts from dominance by a single large basin to a heterogeneous collection of smaller basins. The entropy funnels in Figure~\ref{fig:entropy_funnels} provide a complementary information-theoretic view. Starting from a uniform prior over states, the induced distribution over attractors among those trajectories that have converged exhibits a rapid initial entropy drop, followed by a slower tail. In more concrete terms, most uncertainty about the eventual attractor is resolved within a handful of iterations, but a small residual uncertainty can persist for many steps, especially when many attractors of comparable basin size are present. For $D=3$ and $D=4$, the normalised entropy curves fall close to zero, in line with the near-monopolisation of the state space by a single attractor. For $D=5$ and $D=6$, the residual entropy remains appreciable, reflecting the proliferation of attractors and a more even spread of basin sizes. Thus, entropy funnels capture, in a single scalar time series, both the fast collapse of uncertainty and the dependence of long-term behaviour on basin geometry.\\
~\\
Factoring states into digit multisets yields a first level of symmetry reduction that separates combinatorial and dynamical effects (Figure~\ref{fig:multiset_structure}). On the combinatorial side, both the number of multiset classes and their typical sizes grow quickly with $D$, and the size distributions develop heavier upper tails. This indicates the appearance of rare digit patterns that admit many distinct permutations. On the dynamical side, the distributions of mean distance-to-attractor per multiset remain concentrated on short distances but develop longer tails for $D=5$ and $D=6$. A small subset of multisets is therefore systematically associated with longer transients. 
In other words, as $D$ increases the basin structure is distributed more unevenly across multisets: a few digit patterns dominate large basins, while many others sit in small, long-tailed classes. The analysis of ``easy'' and ``hard'' states in Figure~\ref{fig:easy_vs_hard} moves one step closer to the digit level. By contrasting the fastest and slowest deciles of distance-to-attractor, simple patterns emerge: easy states tend to have larger overall digit spread $g_1$ and, for some digit lengths, higher digit variance. Configurations whose digits are widely dispersed therefore tend, on average, to reach an attractor more quickly. In contrast, the total digit sum carries much less interpretable dynamical information once trivial scaling with $D$ is taken into account. These observations suggest that certain aspects of ``difficulty'' are already visible in low-level digit statistics, even before more refined features are introduced.\\
~\\
Gap space provides a coarser, but highly structured, view of the dynamics. By tracking only the overall spread $g_1$ and the internal gap $g_2$ between the second and third largest digits, Kaprekar's routine induces an empirical first-order Markov chain on a relatively small grid. The resulting flow fields and drift statistics reveal a consistent directional bias: on average, one Kaprekar step tends to increase $g_1$ and decrease $g_2$, pushing states towards regions with a wide overall digit range but a more balanced configuration among the middle digits (Figure~\ref{fig:gap_drift}). For three digits, successor states concentrate near the diagonal $g_1 \approx g_2$, and the coupling between $g_1$ and $g_2$ is strong. As $D$ increases, the stationary distributions shift and the linear relations between $\Delta g_1$ and $g_2$, and between $\Delta g_2$ and $g_1$, weaken in slope. This weakening suggests that, in higher-digit systems, the two gap coordinates no longer suffice to capture the dominant directions of flow, and that additional degrees of freedom in digit space become dynamically relevant. The regression analysis in Figure~\ref{fig:predictability} makes this limitation explicit. For $D=3$, a simple linear model based on $(g_1,g_2,\text{digit sum},\text{digit variance})$ explains a substantial fraction of the variance in distance-to-attractor, with errors on the order of one step. In this regime, low-dimensional digit summaries provide a reasonably informative proxy for dynamical difficulty. For $D=4,5,6$, however, the coefficient of determination $R^2$ drops to a few percent and the learned feature weights shrink and fluctuate in sign. Beyond three digits, the relationship between local digit structure and global convergence time becomes increasingly high-dimensional and entangled, so that linear combinations of a small number of simple features capture only a small part of the behaviour. This breakdown of low-dimensional predictability within the present feature set is consistent with the more fragmented basin geometry and more heterogeneous multiset patterns observed at larger $D$.\\
~\\
Several limitations and extensions follow naturally from these findings. The present study is numerical and restricted to base~10 and digit lengths up to six. From a mathematical perspective, one natural direction is to seek analytic bounds on entropy decay, basin sizes, or gap-space drift as $D$ grows, perhaps by exploiting known structural results on Kaprekar constants and loops in higher bases. Another direction is to refine the coarse-graining schemes: alternative feature sets, nonlinear mappings into gap space, or higher-order Markov projections could be explored to recover part of the predictive power lost at larger $D$. The regression framework could also be extended beyond linear models, for example by investigating whether low-depth decision trees or other simple classifiers find more informative digit combinations without sacrificing interpretability.\\
~\\
More broadly, the picture that emerges is that of a finite information-processing system that rapidly funnels a large set of inputs into a much smaller set of outputs while retaining enough internal structure to defeat simple low-dimensional representations of the dynamics once the state space becomes sufficiently large. Kaprekar dynamics therefore offer a useful model system for studying information funnels and coarse-grained Markov structure in deterministic maps on finite spaces. Similar patterns occur in other settings: in statistical physics, many microscopic configurations are summarised by a few macroscopic variables such as temperature or magnetisation \citep{huang2008statistical, pathria2011theory}; in machine learning, deep networks and information-bottleneck methods restrict information flow through low-dimensional latent or bottleneck layers, yielding compact internal codes that still support accurate predictions \citep{barber2012bayesian, tishby2015deep}; and in decision neuroscience, multiple streams of evidence are modelled as being accumulated into a single decision variable that governs choice \citep{ratcliff2008diffusion,shadlen2001neural,gold2007neural}. Against this background, the Kaprekar map provides a fully tractable model system: the entire state space and all attractors are known explicitly, yet the induced information funnels and coarse-grained dynamics remain non-trivial. This makes the Kaprekar map a natural model system for testing methods that search for informative coarse-grainings or low-dimensional summaries, and for exploring how such techniques might carry over to more biologically motivated dynamical models.

\section*{Supplementary Materials}
No supplementary materials are associated with this article.

\section*{Author contribution}
C.D.D. conceived the study, developed the methodology, implemented the numerical analyses, generated the figures, and wrote and revised the manuscript.

\section*{Funding}
This research was funded by the National Science and Technology Council (NSTC), Taiwan, under grant numbers 112-2410-H-038-027-MY2 and 114-2410-H-038-045.

\section*{Institutional Review Board Statement}
This work is purely computational and does not involve human participants, animals, or identifiable personal data.

\section*{Data Availability Statement}
All data analysed in this study consist of deterministic outputs of the Kaprekar map in base~10 for digit lengths $D=3$--$6$. The MATLAB code used to enumerate the state spaces, compute all summary statistics, and generate the figures is openly available at
\url{https://github.com/ChristophDahl/kaprekar-routine-analysis}
under the MIT licence. The repository includes a master script (\texttt{run\_kaprekar\_all.m}) that reproduces all analyses and figures reported in this manuscript.

\section*{Conflicts of Interest}
The author declares no conflict of interest. The funders had no role in the design of the study; in the collection, analysis, or interpretation of data; in the writing of the manuscript; or in the decision to publish the results.

\bibliographystyle{unsrtnat}  
\bibliography{references_kaprekar}

\end{document}